\input amstex
\documentstyle{amsppt}

\pageheight{21 cm} \pagewidth{16 cm}

\topmatter
\title \break A family of martingales generated by a  process with independent increments
\endtitle

\rightheadtext{Martingales and  independent increment  processes}

\author Josep Llu\'is Sol\'e and Frederic Utzet
\endauthor

\abstract
An explicit procedure to construct a family of martingales generated by a process with independent increments is presented.
The main tools are the polynomials that give the relationship between the moments and cumulants, and
 a set of martingales related to the jumps of the process called Teugels martingales.
\endabstract

\address
Departament de Matem\`atiques,
Facultat de Ci\`encies,
Universitat Aut\'onoma de Barcelona,
08193 Bellaterra (Barcelona),
Spain.
\endaddress

\keywords Process with independent increments, Cumulants, Teugels martingales
\endkeywords

\subjclass 60G51, 60G44
\endsubjclass

\email
jllsole\@mat.uab.cat, utzet\@mat.uab.cat
\endemail

\thanks
This research was supported by grant BFM2006-06247 of
 the Ministerio de Educaci{\'o}n y Ciencia and FEDER
\endthanks

\endtopmatter

\document

\head 1. Introduction
\endhead

In this work we present an explicit procedure to generate a family of martingales from
 a process $X=\{X_t,\, t\ge 0\}$ with independent increments and  continuous in probability.
 We
extend our results exposed in [8], where we dealt with L\'evy processes (independent  and stationary increments);
in that case, the martingales obtained  were of the form $M_t=P(X_t,t)$,  where
$P(x,t)$ is a polynomial in $x$ and $t$, and then they are {\it time--space} harmonic polynomials relative to
$X$.
Here, the martingales constructed are  polynomials on $X_t$ but, in general,  not in $t$.
  Part of the paper is devoted to define the  Teugels martingales of
a process with independent increments; such  martingales, introduced by Nualart and Schoutens [5] for L\'evy
processes, are a building block of the stochastic calculus with that type of  processes.

\head 2. Independent increment  processes and their Teugels martingales
\endhead

Let
$X=\{X_t,\, t\ge 0 \}$ be  a process with independent increments, $X_0=0$,  continuous
in probability and  cadlag; such processes are also called  additive  processes, and we will indistinctly use
    both names. Moreover, assume that $X_t$  is centered, and has moment of all orders.
    It is well known that the law of $X_t$ is infinitely divisible for  all $t\ge 0$.
Let
$\sigma^2_t$  the variance of the gaussian part of $X_t$, and $\nu_t$ its L\'evy measure; for all
these notions we refer to
 Sato [6] or Skorohod [7].

 Denote by $\widetilde \nu$ the (unique) measure on ${\Cal B}((0,\infty)\times \Bbb{R}_0)$ defined by
$$\widetilde \nu((0,t]\times B)=\nu_t(B),\ B\in{\Cal B}(\Bbb{R}_0),$$
 where $\Bbb{R}_0=\Bbb{R}-\{0\}$.
By an standard approximation argument,
we have that for a  measurable function $f:\Bbb{R}_0\to \Bbb{R}$, and for every $t>0,$
$$\iint_{(0,t]\times \Bbb{R}_0}\vert f(x)\vert\, \widetilde \nu(ds,dx)<\infty
 \quad \Longleftrightarrow \quad \int_{\Bbb{R}_0}\vert f(x)\vert\,  \nu_t(dx)<\infty,$$
 and in this case,
$$\iint_{(0,t]\times \Bbb{R}_0} f(x)\, \widetilde \nu(ds,dx)=
  \int_{\Bbb{R}_0} f(x)\,  \nu_t(dx).$$

Note that since  for every $t\ge 0,$ $\nu_t$ is a L\'evy measure, then $\widetilde \nu$ is
$\sigma$--finite. To prove this,  observe that  $\nu_t(\{\vert x\vert >1\})<\infty$, and
$\nu_t((1/(n+1),1/n])<\infty$ and $\nu_t([-1/n,-1/(n+1)))<\infty$, $n\ge 1$. So
there is a numerable partition of $\Bbb{R}_0$ with sets   of finite $\nu_t$ measure, $\forall t>0$. Then,  we can construct a numerable partition of
$(0,\infty)\times \Bbb{R}_0$,  each set with finite $\widetilde \nu$-measure.

Write
$$N(C)=\#\{t:\, (t,\Delta X_t)\in C\},\quad C\in{\Cal B}((0,\infty)\times {\Bbb {R}}_0),$$
 the jump measure of the process, where $\Delta X_t=X_t-X_{t-}$.  It is a Poisson random measure
on $(0,\infty)\times \Bbb{R}_0$ with intensity measure $\widetilde \nu$ (Sato [6, Theorem 19.2]).
 Define the compensated jump measure
 $$d\widetilde N (t,x)=dN(t,x)- d\widetilde\nu(t,x).$$

The process admits the L\'evy--It\^o representation
$$
X_t=G_t+\iint_{(0,t]\times {\Bbb {R}}_0}x\, d\widetilde N(t,x),\eqno{(2)}
$$
where $\{G_t,\, t\ge 0\}$ is a centered continuous Gaussian process with independent increments
and variance
$\Bbb{E}[G^2_t]=\sigma_t^2.$

\bigskip

It is also well known the relationship between the moments of an infinitely divisible law and the {\it moments}
of  its L\'evy measure (see Sato [6,Theorem 25.4]). In our case, as the process has moments of all orders,
for all $t\ge 0,$
$$\int_{\{\vert x\vert >1\}}\vert x\vert\, \nu_t(dx)<\infty\quad \text{and}
\quad\int_{\Bbb{R}_0} \vert x\vert^n\, \nu_t(dx)<\infty,\ \forall n\ge 2.$$
Write

$$
F_2(t)=\sigma_t^2+\int_{\Bbb{R}_0} x^2\, \nu_t(dx) \quad\text{and}\quad
F_n(t)=\int_{\Bbb{R}_0} x^{n}\, \nu_t(dx), \ n\ge 3. \eqno{(1)}
$$
Since $\int_{\{\vert x\vert >1\}}\vert x\vert\, \nu_t(dx)<\infty$ and $\Bbb{E}[X_t]=0$,
the characteristic function of $X_t$ can be written as
$$\phi_t(u)=\exp\big\{- \frac{1}{2}\,\sigma_t^2\,u^2+\int_{\Bbb{R}_0}\big(e^{iux}-1-iux\big)\nu_t(dx)\big\}.$$
It is deduced that  for $n\ge 2$, $F_n(t)$ is the cumulant of order $n$ of $X_t$ (for
$n=2$, $\Bbb{E}[X^2_t]=F_2(t)$). Also $\sigma^2_t$ is continuous and increasing
(Sato [6, Theorem 9.8]).

\proclaim{Proposition 1}
The functions $F_n(t), \ n\ge 2$, are continuous and have finite variation on finite intervals,
and for $n$ even, they are increasing.
\endproclaim

\demo{Proof}

Consider $0<u<t<v$, and write $U=[u,v]$.
From the continuity in probability of $X$,
$$\lim_{s\to t,\, s\in U} X_s^n=X_t^n,\ \hbox{in probability}.$$
Moreover, $\forall s\in U,\ \vert X_s\vert \le \sup_{r\in U}\vert X_r \vert$, and since $X$ is a  martingale, by Doob's inequality,
$$\Bbb{E}\big[\sup_{r\in U}\vert X_r\vert ^n\big]\le C \sup_{r\in U}\Bbb{E}[\vert X_r\vert ^n]\le C\Bbb{E}[\vert X_v\vert ^n]<\infty.$$
So by dominated convergence it follows that the function $t \mapsto \Bbb{E}[X_t^n]$ is continuous.
Since the cumulants are polynomials of the moments, it is deduced the continuity of   all functions
$F_n(t)$.

To prove that $F_n(t)$ has finite variation on finite intervals, consider a partition of $[0,t]$:
 $0<t_0<\cdots<t_k=t.$ Then
 $$\eqalign{
 \sum_{j=1}^{k}\big\vert F_n(t_{j})&-F_n(t_{j-1})\big\vert =
\sum_{j=1}^{k} \big\vert\iint_{(t_{j-1},t_j]\times\Bbb{R}_0} x^n\widetilde\nu(ds,dx)\big\vert\cr
&\le \sum_{j=1}^{k} \iint_{(t_{j-1},t_j]\times\Bbb{R}_0} \vert x\vert ^n\,\widetilde\nu(ds,dx)
 =  \iint_{(0,t]\times\Bbb{R}_0} \vert x\vert^n\,\widetilde\nu(ds,dx)<\infty.\quad \blacksquare}$$
 \enddemo

\bigskip

 Consider  the {\it variations} of the process $X$ (see Meyer [4]):
$$\eqalign{
X^{(1)}_t&=X_t, \cr
X ^{(2)}_t&=[X,X]_t=\sigma^2_t+\sum _{0<s\le t} \big(\Delta X_s \big)^2\cr
X ^{(n)}_t&=\sum _{0<s\le t} \big(\Delta X_s \big)^n, \  n \ge 3.\cr}$$
By Kyprianou [2, Theorem 2.7], for $n\ge 3$ (the case $n=2$ is similar), the characteristic function  of $X^{(n)}$
is
$$\exp\big\{\iint_{(0,t]\times {\Bbb{R}}_0}\big(e^{i ux^n}-1\big)\, \widetilde \nu(ds,dx)\,\big\}
=\exp\big\{\int_{\Bbb{R}_0}\big(e^{i ux}-1\big)\, \nu_t^{(n)}(dx)\,\big\},
$$
where $\nu_t^{(n)}$ is the measure image of $\nu_t$ by the function  $x \mapsto x^n$, which is  a L\'evy measure.
So $X^{(n)}$  has independent increments. Also by Kyprianou [2, Theorem 2.7], for $n\ge 2$,
$$\Bbb{E}[X^{(n)}_t]=F_n(t)\quad \hbox{and}\quad \Bbb{E}\big[\big(X^{(n)}_t\big)^2\big]=F_{2n}(t)+\big(F_n(t)\big)^2.$$
Therefore, combining the independence of the increments  and the continuity of $F_n(t)$,
it is  deduced that $X^{(n)}$ is continuous in probability.

By Proposition 1, $F_n(t)$ has finite variation on finite intervals. Hence, the process $$X_t^{(n)}=F_n(t)+\big(X_t^{(n)}
-F_n(t)\big)$$ is a semimartingale.

\bigskip

 \bigskip

 The Teugels martingales introduced by  Nualart and Schoutens [5] for L\'evy processes can be
 extended to additive processes. In the same way as in  [5], these  martingales are obtained
 centering the processes $X^{(n)}$:
$$\eqalign{
Y^{(1)}_t&=X_t,\cr
Y^{(n)}_t&=X^{(n)}-F_n(t),\ n\ge 2,\cr}$$
They are square integrable martingales with optional quadratic covariation
$$[Y^{(n)},Y^{(m)}]_t=X^{(n+m)},$$
and, since $F_{2n}(t)$ is increasing, the predictable quadratic variation of $Y^{(n)}$ is
$$\langle Y^{(n)}\rangle_t=F_{2n}(t).$$

\head 3. The  polynomials of cumulants
\endhead

 The formal expression
$$\exp\bigg\{\sum_{n=1}^\infty\kappa_n\frac{u^n}{n!}\bigg\}=\sum_{n=0}^\infty\mu_n\frac{u^n}{n!}. \eqno(3)$$
relates the sequences of numbers $\{\kappa_n,\ n\ge 1\}$ and $\{\mu_n,\ n\ge 0\}$.
When we consider a random variable $Z$  with moment generating function in some open interval containing 0,
then  both series converge in a neighborhood of 0, and (3) is the relationship between
 the moment generating
function, $\psi(u)=\Bbb{E}[e^{uZ}]$, and the cumulant generating function, $\log \psi (u)$. Moreover,
$\mu_n$ (respectively,  $\kappa_n$)
 is the moment  (respectively, the cumulant) of order $n$ of $Z$, and the well known
 relationship between moments and cumulants can be deduced from (3). The first three are
$$\eqalign{\mu_1&=\kappa_1,\cr
  \mu_2&=\kappa_1^2+\kappa_2,\cr
  \mu_3&=\kappa_1^3+3\kappa_1\kappa_2+\kappa_3, \dots\cr}$$
If the random variable $Z$ has only finite moments up to order $n$,  the corresponding relationship is true up to this order.

There is a general explicit expression of the moments in terms of cumulants in Kendall and Stuart [1], or
formulas  involving the partitions of a set,  see McCullagh [3]. In general, $\mu_n$ is a polynomial
of $\kappa_1,\dots,\kappa_n$, called Kendall polynomial.
Denote by $\Gamma_n(x_1,\dots,x_n)$, $n\ge 1$, this  polynomial, that is, we have
$$\mu_n=\Gamma_n(\kappa_1,\dots,\kappa_n).$$
Also write $\Gamma_0=1$. These polynomials enjoy very interesting properties, as
the recurrence formula that  follows from   Stanley [9, Proposition 5.1.7]:

$$\Gamma_{n+1}(x_1,\dots,x_{n+1})=\sum_{j=0}^n{n\choose j}\Gamma_j(x_1,\dots,x_j)\,x_{n+1-j}.\eqno(4)$$
We  also have
$$
\frac{\partial \Gamma _n (x_1,\dots,x_n)}{\partial x_j}={n\choose j} \Gamma_{n-j}(x_1,\dots,x_{n-j}),\quad j=1,\dots,n.\eqno(5)
$$
 Computing the Taylor expansion of $\Gamma _n (x_1+y,x_2,\dots,x_n)$   at $y=0$, we get the following expression that we will
 need later:
$$
\Gamma_n(x_1+y,x_2,\dots,x_n)=\sum_{j=0}^n {n\choose j}\Gamma_{n-j}(x_1,\dots,x_{n-j})\,y^{j}.\eqno(6)
$$
Interchanging the role of $x_1$ and $y$, and evaluating the function at 0 we obtain
$$
\Gamma_n(x_1,x_2,\dots,x_n)=\sum_{j=0}^n {n\choose j}\Gamma_{n-j}(0,x_2,\dots,x_{n-j})\,x_1^{j}. \eqno(7)
$$

\head 4. A family of martingales  relative to the additive process
\endhead

The  main result of the paper is the following Theorem:
\proclaim{Theorem 1}
Let $X$ be a centered additive process with finite moments of all orders. Then  the process
$$M^{(n)}_t=\Gamma_n\Big(X_t,- F_2(t),\dots, -F_n(t) \Big)$$
is a martingale.
\endproclaim

\demo{Proof}

Let $n\ge 2$. We apply the multidimensional It\^o   formula  to the semimartingales $X_t,\, F_2(t),\dots, F_n(t)$.
By Proposition 1,  the functions $F_2(t),\dots, F_n(t)$ and $\sigma_t^2$ are continuous and of finite variation.
From (5) and the fact that $[X,X]^c_t=\sigma_t^2$ and $[F_j,F_j]^c_t=0$, we have
$$\eqalign{
M^{(n)}_t& =n\int_0^tM^{(n-1)}_{s-}\, dX_s-\sum_{j=2}^n
{n\choose j}\int_0^t  M^{(n-j)}_s dF_j(s)+\frac{1}{2} n(n-1)\int_0^t M^{(n-2)}_s\, d(\sigma^2_s) \cr
& \qquad
 + \sum_{0<s\le t}\Big( \Gamma_n\big(X_{s-}+\Delta X_s,- F_2(s),\dots, -F_n (s) \big)
 - \Gamma_n\big(X_{s-},- F_2(s),\dots, - F_n(s) \big)\cr
& \qquad -
n \Delta X_s \,\Gamma_{n-1}\big(X_{s-},-F_2(s),\dots, -F_n(s) \big)\Big)\cr}
$$
Applying (6),
$$\Gamma_n\big(X_{s-}+\Delta X_s,-F_2(s),\dots, -F_n(s)\big)=
\sum_{j=0}^n{n\choose j}M^{(n-j)}_{s-}\big(\Delta X_s\big)^j.
$$
Then, the jumps part given in the expression of $M^{(n)}_t$
is
$$\eqalign{
\sum_{0<s\le t} \sum_{j=2}^n {n\choose j}&M^{(n-j)}_{s-}\big(\Delta X_s\big)^j=
\sum_{j=2}^n{n\choose j}\int_0^t M^{(n-j)}_{s-} \,d X_s^{(j)}-{n\choose 2}\int_0^t M^{(n-2)}_s\, d(\sigma^2_s)\cr
&=\sum_{j=2}^n{n\choose j}\int_0^t M^{(n-j)}_{s-} \,
  d\big( Y_s^{(j)}+ F_j(s)\big)-{n\choose 2}\int_0^t M^{(n-2)}_s\, d(\sigma^2_s). \cr}
$$
Therefore,
$$
M^{(n)}_t=\sum_{j=1}^n{n\choose j}\int_0^tM^{(n-j)}_{s-}\, dY^{(j)}_s.\eqno{(8)}$$
Moreover, $\big(M^{(k)}_t\big)^2$ is a polynomial in $X_t,\, F_2(t),\dots, F_k(t)$. Taking expectations and using
the  relationship between moments and cumulants, and that the cumulants of $X_t$ are $F_n(t)$, $n\ge 2$,
we obtain that
$$E\big[\big(M^{(k)}_t\big)^2\big]=P(F_2(t),\dots,F_{2k}(t)),$$
for a suitable polynomial $P$.
Then, for every $t\ge 0$, we have
$$\Bbb{E}\big[\int_0^t \big(M^{(k)}_{s-}\big)^2\, d\langle Y^{(j)}\rangle_s\big]=\int_0^t\Bbb{E}\Big[\big(M^{(k)}_{s-}\big)^2\Big] d F_{2j}(s)=
\int_0^t P(F_2(s),\dots,F_{2k}(s))\, d F_{2j}(s)
<\infty,$$
So all the stochastic integrals on the right hand side of (8) are martingales.\quad $\blacksquare$

\enddemo

\proclaim{Remark 1} It is worth to note that it follows from the preceding Theorem that the function
$$g_n(x,t)=\Gamma_n\Big(x,- F_2(t),\dots, -F_n(t) \Big)$$
is a time--space harmonic function with respect to $X_t$.  By (7)
$$g_n(x,t)=\sum_{j=0}^n\Gamma_{n-j}(0,-F_2(t),\dots, F_n(t))x^j.$$
In general, $g_n(x,t)$  is a polynomial in $x$.  If  $F_n(t),\ n\ge 2,$ are polynomials
in $t$, then $g_n(x,t)$ is a time--space harmonic polynomial; this happens for all L\'evy processes with moments
of all orders,
and for  some  additive process; see the example below.
\endproclaim

\noindent {\bf Example.} Let $\Lambda(t):\Bbb{R}_+\longrightarrow \Bbb{R}_+$ be a continuous increasing function,
and $J$ be a Poisson random measure  on $\Bbb{R}_+$ with intensity measure
 $\mu(A)=\int_A\Lambda(dt),\ A\in \Cal{B}(\Bbb{R}_+)$. Then the process $X=\{X_t,\, t\ge 0\}$ defined pathwise,
 $$X_t(\omega)=\int_0^t J(ds,\omega)-\Lambda(t)$$
 is an additive process; it is  a Cox process with deterministic hazard function $\Lambda(t)$.
 From the characteristic function of $X_t$ we deduce that the L\'evy measure is
 $$\nu_t(dx)=\Lambda(t)\delta_1(dx),$$
 where $\delta_1$ is a Dirac delta measure concentrated in the point 1. Hence,
 $$F_n(t)=\Lambda(t), \ n\ge 2.$$
 Note that the conditions that we have assumed  on $\Lambda$ are necessary to obtain an additive process, but it is not necessary
 (though not  very restrictive) to assume that $\Lambda$ is absolutely continuous with respect to the Lebesgue measure.

 The function defined in Remark 1 is
 $$g_n(x,t)=\Gamma_n\big(x,-\Lambda(t),\dots, -\Lambda(t)\big).$$
Hence,  when $\Lambda(t)$ is a polynomial,   $g_n(x,t)$ is a time--space harmonic polynomial.

Denote by $\overline C _n(x,t)$ the Charlier polynomial with leading coefficient equal to 1, then (see [8])
$$g_n(x,t)=\sum_{j=1}^n\lambda^{(n)}_j \overline C_j(x,\Lambda(t)),$$
where
$\lambda^{(n)}_1=1$ and
$$\lambda^{(n)}_{k+1}=\sum_{j=k}^{n-1}{n\choose j}\lambda^{(j)}_k,\ k=1,\dots, n-1.$$


\Refs



\ref
\no1
\by M. Kendall  and  A. Stuart
\book The Advanced Theory of Statistics, Vol. 1, 4th edition
\publ MacMillan Publishing CO., Inc.
\publaddr New York
\yr 1977
\endref

\ref
\no2
\by A. E. Kyprianou
\book Introductory Lectures on Fluctuations of L\'evy Processes with Applications
\publ Springer
\publaddr Berlin
\yr 2006
\endref



\ref
\no3
\by P. McCullagh
\book Tensor Methods in Statistics
\publ Chapman and Hall
\publaddr London
\yr 1987
\endref

\ref
 \no4
 \by P. A. Meyer
 \paper Un cours sur les integrales stochastiques
 \inbook S\'eminaire de Probabilit\'es X
\publ Springer
 \publaddr New York
 \yr 1976 \pages245--400
  \lang  French
\endref




\ref
 \no5
  \by D. Nualart and  W. Schoutens
  \paper Chaotic and predictable representation for L\'{e}vy processes
\jour Stochastic Process. Appl.
\vol 90
\yr 2000
\pages 109--122
\endref




\ref
\no6
\by K. Sato
\book L\'{e}vy Processes and Infinitely Divisible Distributions
\publ Cambridge University Press
\publaddr Cambridge
\yr 1999
\endref

\ref
\no7
\by A. V. Skorohod
\book Random Processes with Independent Increments
\publ Kluwer Academic Publ.
\publaddr Dordrecht, Boston, London
\yr 1986
\endref




\ref
 \no 8
  \by J. L. Sol\'e and  F. Utzet
  \paper Time--space harmonic polynomials relative to a L\'evy process
\jour Bernoulli
\yr 2007
\endref







\ref
\no9
\by R. P. Stanley
\book Enumerative Combinatorics, Vol. 2
\publ Cambridge University Press
\publaddr Cambridge
\yr 1999
\endref

\endRefs
\enddocument